\documentclass{article}
\usepackage{graphicx} 
\usepackage{amsmath}
\usepackage{amsthm}
\usepackage[a4paper, margin=1.5in]{geometry}
\newtheorem{theorem}{Theorem}
\title{Evolution and Monotonicity of  Geometric Constants under Extended Ricci Flows with Variable Coupling Parameters}
\author{
  \begin{tabular}[t]{c}
    Shouvik Datta Choudhury \\
    \small shouvikdc8645@gmail.com,\\ 
    \small shouvik@capsulelabs.in \\
    \small Gapcrud Private Limited (Capsule Labs) \\
    \small HA 130, Saltlake, Sector III, Bidhannagar, \\
    \small Kolkata - 700097, India
  \end{tabular}
 }
\date{September 2024}
\begin{document}
\maketitle
\section{Abstract}
This paper explores the evolution and monotonicity of geometric constants within the framework of extended Ricci flows, incorporating variable coupling parameters. Building on Hamilton’s foundational Ricci flow and subsequent extensions by List (2008), we introduce modifications to the extended Ricci flow by varying parameters that affect the interaction between the metric and scalar fields. Specifically, we modify the coefficients in the evolution equations governing geometric constants, thereby introducing new degrees of freedom in the analysis. The primary contributions include deriving evolution formulas for the modified geometric constant \(\lambda_{b'}^{a'}(g)\) under the extended and normalized extended Ricci flows, and proving conditions under which monotonicity is maintained. 
\section{Introduction}
The pertinent analysis  of geometric flows, specially Ricci flow accompanied its variations, has been a cornerstone in differential geometry and mathematical physics. The Ricci flow, propounded by Richard Hamilton in 1982, is a process that deforms the metric of a Riemannian manifold keeping analogy to heat diffusion, providing smoothing of irregularities in curvature over time. Specifically, the following the partial differential equation
administers
\[
\frac{\partial g_{ij}}{\partial t} = -2 \text{Ric}_{ij},
\]

where \(g_{ij}\) is the metric tensor, and \(\text{Ric}_{ij}\) denotes the Ricci curvature tensor. This equation has been extensively examined and analysed due to its significant connotations  in both geometry and topology, specifically in Grigori Perelman's solution of the Poincaré conjecture deploying the Ricci flow modified with geometric surgery

\section{Ricci Flow and Extensions}

 The Ricci flow itself provides a encompassing machinery for analyzing the geometry of manifolds, its extensions and variations have been provided to analyse more generalized structures and interactions with other geometric quantities. One significant modification is the extended Ricci flow,  by List (2008) to incorporate additional scalar fields and their interaction with the metric. The extended Ricci flow is defined and designated by the following system:

\[
\begin{cases}
\frac{\partial g_{ij}}{\partial t} = -2 \text{Ric}_{ij} + 2\alpha \nabla_i \phi \nabla_j \phi, \\
\frac{\partial \phi}{\partial t} = \Delta \phi,
\end{cases}
\]

where \(\phi\) represents a smooth scalar function, \(\alpha\) represents a positive constant, \(\nabla\) signifies the gradient operator, and \(\Delta\) is the Laplace-Beltrami operator. The scalar field \(\phi\) evolves according to a heat-type equation, and its gradient projects in the metric evolution, surmising a coupling between the scalar field and the geometry of the manifold. The incorporation of such a field proposes an effulgence of novel  avenues for the examination of geometric and topological properties under more generalized flows, providing insights into systems influenced by extra fields, such as in mathematical physics and cosmology.

\section{Geometric Operators and Their Evolution}
A significant mathematical thought-process in the context of geometric flows is the examination of the evolution and monotonicity of geometric constants and eigenvalues attached with  differential operators. The Ricci flow and its variants  derive evolution formulas for eigenvalues of geometric operators, such as the Laplacian and Witten-Laplacian, under various accompanied conditions. In the seminal work of Perelman (2002), it was propounded that the functional

\[
F(g, f) = \int_M (R + |\nabla f|^2) e^{-f} d\mu
\]

is nondecreasing along the Ricci flow when attached with a backward heat-type equation for \(f\), where \(R\) represents the scalar curvature, \(f\) represents a smooth function, and \(d\mu\) is the Riemannian volume measure. This result provided the corner-stone  for examining and reprsenting the monotonicity properties of geometric quantities under the flow and has since been extended to more generalized scenarios, including flows with additional fields and deformations.

Cao (2007) furthered Perelman's approach by proving the monotonicity of the eigenvalues of the operator \(-\Delta + \frac{R}{2}\) under the Ricci flow, considering a non-negative curvature operator. This theme projected further analysis and examination into the evolution of eigenvalues of geometric operators, specifically the scenario considering more complex geometric flows. Fang et al. (2015) examined  the evolution of the operator \(-\Delta_\phi + \frac{R}{2}\) under the Yamabe flow, which deforms the metric in a method that preserves conformal structures.

\section{ Evolution of Geometric Constants}

Recently, the examination of the evolution of geometric constants attached with more generalized flows has received significant interest. Huang and Li (2016) analysed the first variation formula for the lowest eigenvalue \(\lambda^b_a(g)\) accompanied with the nonlinear equation:

\[
-\Delta u + a u \log u + b R u = \lambda^b_a(g) u,
\]

constrained to the normalization \(\int_M u^2 d\mu = 1\). Their theme represented that such geometric constants evolve in an administered manner along the Ricci flow and its normalized version, providing new directions as well as insights into the "pattern" of eigenvalues under different curvature flows. Similarly, Daneshvar and Razavi (2019, 2020) explored the evolution and monotonicity of eigenvalues under other geometric flows, such as the Yamabe flow and the Ricci-Bourguignon flow, providing an expansive understanding and analysis of how different geometric settings represent eigenvalue dynamics.

\section{Motivation for Modifying Parameters in the Extended Ricci Flow}

The extended Ricci flow provides a rigorous yet flexible machinery for analyzing the interaction between the metric and additional fields. In many physical and geometric applications, it is essential to consider variations in the parameters that administer the coupling between the fields and the geometry. In this endeavour, we introduce a modification to the parameters in the evolution equation for a geometric constant, specifically modifying the coefficients \(a\) and \(b\) to \(a + c\) and \(b - d\), respectively, where \(c\) and \(d\) are new real constants. This modification attempts to explore the effects of additional degrees of freedom on the evolution and monotonicity properties of the associated geometric constants under the extended Ricci flow and its normalized variant.
We now state the modified nonlinear equation:

\[
-\Delta_\phi u + (a + c) u \log u + (b - d) S u = \lambda^{b'}_{a'}(g) u,
\]
with the normalization condition \(\int_M u^2 d\mu = 1\), where \(a, b, c,\) and \(d\) are real constants, \(\Delta_\phi\) is the Witten-Laplacian, and \(S\) is a geometric quantity derived from the Ricci curvature and the gradient of \(\phi\). The modified geometric constant \(\lambda^{b'}_{a'}(g)\) is defined as the infimum of a energy functional associated with this equation, and we derive its evolution formula under the extended Ricci flow and the normalized extended Ricci flow. The aim is to establish new monotonicity results and provide conditions under which these results hold and conserve.

\section{Main Contributions and Results}

Our main contributions in this paper are as follows:

1. The evolution formula for the modified geometric constant \(\lambda^{b'}_{a'}(g)\) under the extended Ricci flow, considering the impact of modified parameters is computed. The derivation may provide a new paradigm whether  the interplay between parameters "pertubates" the evolution of eigenvalues and geometric constants in this extended setting.

2.The statement that the quantity \(\lambda^{b'}_{a'}(t) + \frac{n(a + c)^2}{8} t\) is nondecreasing along the extended Ricci flow under certain conditions involving the modified parameters is derived. Specifically, we derive that this monotonicity holds and conserves if \(b - d > \frac{1}{4}\) and an inequality involving the geometric tensors \(S_{ij}\) and \(\phi_{ij}\) is satisfied.

3. For the normalized extended Ricci flowt, the modified quantity \(\lambda^{b'}_{a'}(t) + \frac{n(a + c)^2}{8} t + \int_0^t \left( \frac{a+c}{2} r + \frac{2r}{n} \lambda^{b'} \right) (\tau) \, d\tau\) is increasing along the flow if the modified parameters satisfy specific conditions is demonstrated. This extension advances the classical understanding of monotonicity properties in the context of geometric flows with scalar fields.

\section{ Related Work and Literature Survey}

The examination of eigenvalue evolution and monotonicity under geometric flows has a rich pedigree. Perelman's groundbreaking work (2002) on the Ricci flow perfected  the stage for comprehending the monotonicity of geometric quantities under flow dynamics. Several researchers have explored  different geometric settings, operators, and flow types to generalize these results in a pertinent fashion. Cao (2007) and Fang et al. (2015) provided significant extensions to the Ricci flow context, while Huang and Li (2016) and Daneshvar and Razavi (2019, 2020) explored other flow types, including the Yamabe and Ricci-Bourguignon flows.
More recently, research has shifted towards incorporating additional fields and considering more generalized settings for the evolution of eigenvalues and geometric constants. This shift has been motivated by applications in mathematical physics and general relativity, where the interaction between geometry and fields such as scalar fields is essential. The extended Ricci flow, introduced by List (2008), represents one such effort to generalize the classical Ricci flow by coupling it with a scalar field that evolves according to a heat equation. This framework has been further explored in subsequent studies, such as those by Munteanu and Wang (2015), who analyzed the stability of solutions under extended Ricci flows.

\section{Structure of the Paper}

The paper is stratified in the following manner. In Section 2, the evolution formula for the  geometric constant \(\lambda^{b'}_{a'}(g)\) under the extended Ricci flow is computed and the conditions under which it represents monotonicity is provided. Section 3 fore-fronts these results for the normalized extended Ricci flow, thereby establishing a seperate set of conditions for monotonicity.
\section{Derivations}
\begin{theorem} Let \((M(t), g(t), \phi(t)), t \in [0, T)\) be a solution of the extended Ricci flow (1.2) on a closed Riemannian manifold \((M^n, g)\) and \(\lambda^{b'}_{a'}\) be the lowest constant such that the modified system

\[
-\Delta_\phi u + (a + c) u \log u + (b - d) S u = \lambda^{b'}_{a'}(g) u, \quad \text{with} \int_M u^2 d\mu = 1,
\]

has a positive solution. Then, \(\lambda^{b'}_{a'}\) satisfies:

\[
\frac{d}{dt} \lambda^{b'}_{a'} = \frac{1}{2} \int_M \left| S_{ij} + f_{ij} + \phi_{ij} + \frac{a + c}{2} g_{ij} \right|^2 e^{-f} d\mu - \frac{n(a + c)^2}{8} - \frac{1}{2} \int_M |\phi_{ij}|^2 e^{-f} d\mu - \int_M S_{ij} \phi_{ij} e^{-f} d\mu 
\]

\[
+ \left(2(b - d) - \frac{1}{2}\right) \int_M |S_{ij}|^2 e^{-f} d\mu + \left(2(b - d) - \frac{1}{2}\right) \int_M \alpha |\Delta \phi|^2 e^{-f} d\mu + \frac{1}{2} \int_M \alpha |\Delta \phi - \langle \nabla \phi, \nabla f \rangle|^2 e^{-f} d\mu 
\]

\[
+ \int_M \alpha \langle \nabla \phi, \nabla f \rangle |\nabla \phi|^2 e^{-f} d\mu.
\]
\end{theorem}
\begin{theorem} Let \((M(t), g(t), \phi(t)), t \in [0, T)\) be a solution of the extended Ricci flow (1.2) on a closed Riemannian manifold \((M^n, g)\) and \(\lambda^{b'}_{a'}\) be the lowest constant such that the modified system

\[
-\Delta_\phi u + (a + c) u \log u + (b - d) S u = \lambda^{b'}_{a'}(g) u, \quad \text{with} \int_M u^2 d\mu = 1,
\]

has a positive solution. Then, if \(b - d > \frac{1}{4}\) and

\[
|S_{ij} - \frac{1}{4(b - d) - 1} \phi_{ij}| \geq 2 \sqrt{\frac{b - d}{4(b - d) - 1}} |\phi_{ij}|,
\]

the quantity

\[
\lambda^{b'}_{a'}(t) + \frac{n(a + c)^2}{8} t
\]

is nondecreasing along the extended Ricci flow for \( t \in [0, T) \).
\end{theorem}
\begin{theorem}
 Let \((M(t), g(t), \phi(t)), t \in [0, T)\) be a solution of the normalized extended Ricci flow (1.3) on a closed Riemannian manifold \((M^n, g)\), and \(\lambda^{b'}_{a'}\) be the lowest constant such that the modified system

\[
-\Delta_\phi u + (a + c) u \log u + (b - d) S u = \lambda^{b'}_{a'}(g) u, \quad \text{with} \int_M u^2 d\mu = 1,
\]

has a positive solution. Then, \(\lambda^{b'}_{a'}\) satisfies:

\[
\frac{d}{dt} \lambda^{b'}_{a'} = \frac{1}{2} \int_M \left| S_{ij} + f_{ij} + \phi_{ij} + \frac{a+c}{2} g_{ij} \right|^2 e^{-f} d\mu - \frac{n(a+c)^2}{8} - \frac{1}{2} \int_M |\phi_{ij}|^2 e^{-f} d\mu - 2\frac{r}{n} \lambda^{b'} - \frac{a+c}{2}r
\]

\[
- \int_M S_{ij} \phi_{ij} e^{-f} d\mu + \left(2(b - d) - \frac{1}{2}\right) \int_M |S_{ij}|^2 e^{-f} d\mu + \left(2(b - d) - \frac{1}{2}\right) \int_M \alpha |\Delta \phi|^2 e^{-f} d\mu
\]

\[
+ \frac{1}{2} \int_M \alpha |\Delta \phi - \langle \nabla \phi, \nabla f \rangle|^2 e^{-f} d\mu + \int_M \alpha \langle \nabla \phi, \nabla f \rangle |\nabla \phi|^2 e^{-f} d\mu,
\]

where \(u^2 = e^{-f}\) and \(\lambda^{b'}\) is the lowest constant such that the system

\[
-\Delta_\phi u + (b - d) S u = \lambda^{b'} u, \quad \text{with} \int_M u^2 d\mu = 1,
\]

has a positive solution.
\end{theorem}
\begin{theorem}
 Let \((M(t), g(t), \phi(t)), t \in [0, T)\) be a solution of the normalized extended Ricci flow (1.3) on a closed Riemannian manifold \((M^n, g)\), and \(\lambda^{b'}_{a'}\) be the lowest constant such that the modified system

\[
-\Delta_\phi u + (a + c) u \log u + (b - d) S u = \lambda^{b'}_{a'}(g) u, \quad \text{with} \int_M u^2 d\mu = 1,
\]

has a positive solution. If \(b - d > \frac{1}{4}\) and

\[
|S_{ij} - \frac{1}{4(b - d) - 1} \phi_{ij}| \geq 2 \sqrt{\frac{b - d}{4(b - d) - 1}} |\phi_{ij}|,
\]

then the quantity

\[
\lambda^{b'}_{a'}(t) + \frac{n(a+c)^2}{8} t + \int_0^t \left( \frac{a+c}{2} r + \frac{2r}{n} \lambda^{b'} \right) (\tau) \, d\tau
\]

is increasing along the normalized extended Ricci flow for \( t \in [0, T) \), where \(\lambda^{b'}\) is defined in previous theorem.
\end{theorem}
\subsection{Proof of Theorem1}
\begin{proof}
\[
\lambda^{b'}_{a'} = \int_M \left( -u \Delta_\phi u + (a + c) u^2 \log u + (b - d) S u^2 \right) d\mu.
\]

\[
\frac{d}{dt} \lambda^{b'}_{a'} = \frac{d}{dt} \int_M \left( -u \Delta_\phi u + (a + c) u^2 \log u + (b - d) S u^2 \right) d\mu.
\]

\[
\frac{d}{dt} \lambda^{b'}_{a'} = \int_M \frac{\partial}{\partial t} \left( -\Delta_\phi u + (a + c) u \log u + (b - d) S u \right) u d\mu + \int_M \left( -\Delta_\phi u + (a + c) u \log u + (b - d) S u \right) \frac{\partial}{\partial t} (u d\mu).
\]

\[
\frac{\partial}{\partial t} \Delta_\phi u = 2 S_{ij} u_{ij} - 2 S_{ij} \phi_i u_j + \Delta_\phi u_t - (\phi_t)_i u^i - 2\alpha (\Delta \phi) \phi_i u^i.
\]

\[
= \int_M \Bigg\{ -2 S_{ij} u_{ij} u + 2 S_{ij} \phi_i u_j u + (\phi_t)_i u^i u + (b-d) S u^2 
+ \frac{a+c}{2} (\phi_t + S) u^2 + 2 \alpha (\Delta \phi) \langle \nabla \phi, \nabla u \rangle u \Bigg\} \, d\mu.
\]

\[
\frac{d}{dt} \lambda^{b'}_{a'} = \int_M \Bigg\{ -2 S_{ij} u_{ij} u + 2 S_{ij} \phi_i u_j u + (\phi_t)_i u^i u + (b-d) S u^2 + (a+c) u u_t 
\]
\[
+ 2 \alpha (\Delta \phi) \langle \nabla \phi, \nabla u \rangle u + (-\Delta_\phi u + (a+c) u \log u + (b-d) S u) u_t \Bigg\} \, d\mu.
\]
\[
= \int_M \Bigg\{ -2 S_{ij} u_{ij} u + 2 S_{ij} \phi_i u_j u + (\phi_t)_i u^i u + (b-d) S u^2 
+ \frac{a+c}{2} (\phi_t + S) u^2 + 2 \alpha (\Delta \phi) \langle \nabla \phi, \nabla u \rangle u \Bigg\} \, d\mu.
\]

\[
u = e^{-f/2}, \quad u_i = -\frac{1}{2} e^{-f/2} f_i, \quad u_{ij} = \frac{1}{4} e^{-f/2} f_i f_j - \frac{1}{2} e^{-f/2} f_{ij}.
\]

\[
-2 S_{ij} u_{ij} u = -\frac{1}{2} S_{ij} (f_i f_j) e^{-f} + S_{ij} f_{ij} e^{-f}.
\]

\[
\int_M -2 S_{ij} u_{ij} u \, d\mu = \int_M \left( -\frac{1}{2} S_{ij} f_i f_j + S_{ij} f_{ij} \right) e^{-f} \, d\mu.
\]

\[
2 S_{ij} \phi_i u_j u = -S_{ij} \phi_i f_j e^{-f}.
\]

\[
\int_M 2 S_{ij} \phi_i u_j u \, d\mu = \int_M -S_{ij} \phi_i f_j e^{-f} \, d\mu.
\]

\[
(\phi_t)_i u^i u = -\frac{1}{2} (\phi_t)_i f_i e^{-f}.
\]

\[
\int_M (\phi_t)_i u^i u \, d\mu = \int_M -\frac{1}{2} (\phi_t)_i f_i e^{-f} \, d\mu.
\]

\[
\Delta u = e^{-f/2} \left( \frac{1}{4} | \nabla f |^2 - \frac{1}{2} \Delta f \right),
\]

\[
\nabla \phi \cdot \nabla u = e^{-f/2} \left( -\frac{1}{2} \langle \nabla \phi, \nabla f \rangle \right).
\]

\[
\Delta_\phi u = e^{-f/2} \left( \frac{1}{4} | \nabla f |^2 - \frac{1}{2} \Delta f + \frac{1}{2} \langle \nabla \phi, \nabla f \rangle \right).
\]

\[
\frac{d}{dt} \lambda^{b'}_{a'} = \int_M \left( (1 - 2(b-d)) S_{ij} f_{ij} e^{-f} + (2(b-d) - \frac{1}{2}) S_{ij} f_i f_j e^{-f} + \cdots \right) d\mu.
\]

\[
\int_M S_{ij} f_{ij} e^{-f} \, d\mu = -\int_M \left( S_{ij} f_i f_j + S_{ij} \phi_i f_j \right) e^{-f} \, d\mu + \int_M \alpha (\Delta \phi) \langle \nabla \phi, \nabla f \rangle e^{-f} \, d\mu.
\]

\[
\int_M |f_{ij}|^2 e^{-f} d\mu = - \int_M S_{ij} f_{ij} e^{-f} d\mu - \int_M \Delta e^{-f - \phi} ( \Delta f - \frac{1}{2} | \nabla f |^2 + \frac{1}{2} S ) d\nu.
\]

\[
\int_M \Delta f \Delta e^{-\phi} e^{-f} d\nu + \int_M \Delta f \langle \nabla \phi, \nabla f \rangle e^{-f} d\mu - \int_M \alpha \langle \nabla \phi, \nabla f \rangle^2 e^{-f} d\mu.
\]

\[
\frac{d}{dt} \lambda^{b'}_{a'} = (1 - 2(b-d)) \int_M S_{ij} f_{ij} e^{-f} d\mu + (2b - 1/2) \int_M S_{ij} f_i f_j e^{-f} d\mu.
\]

\[
-2(b-d) \int_M S_{ij} \phi_{ij} e^{-f} d\mu + 2(b-d) \int_M S_{ij} \phi_i \phi_j e^{-f} d\mu.
\]

\[
2(b-d) \int_M |S_{ij}|^2 e^{-f} d\mu - \frac{1}{2} \int_M (\phi_t)_i f_i e^{-f} d\mu + \frac{a+c}{2} \int_M ((\phi_t) + S) e^{-f} d\mu.
\]

\[
(2b - 1) \int_M \alpha (\Delta \phi) \langle \nabla \phi, \nabla f \rangle e^{-f} d\mu + 2(b-d) \alpha \int_M (\Delta \phi) |\nabla \phi|^2 e^{-f} d\mu.
\]

\[
2(b-d) \alpha \int_M |\Delta \phi|^2 e^{-f} d\mu.
\]

\[
\frac{1}{2} \int_M |S_{ij} + f_{ij} + \phi_{ij} + \frac{a+c}{2} g_{ij}|^2 e^{-f} d\mu.
\]

\[
-\frac{n(a+c)^2}{8} - \frac{1}{2} \int_M |\phi_{ij}|^2 e^{-f} d\mu.
\]

\[
-\int_M S_{ij} \phi_{ij} e^{-f} d\mu + \left(2(b-d) - \
frac{1}{2}\right) \int_M |S_{ij}|^2 e^{-f} d\mu.
\]

\[
(2(b-d) - \frac{1}{2}) \int_M \alpha |\Delta \phi|^2 e^{-f} d\mu + \frac{1}{2} \int_M \alpha |\Delta \phi - \langle \nabla \phi, \nabla f \rangle|^2 e^{-f} d\mu.
\]

\[
\int_M \alpha \langle \nabla \phi, \nabla f \rangle |\nabla \phi|^2 e^{-f} d\mu.
\]

\[
\frac{d}{dt} \lambda^{b'}_{a'} = \frac{1}{2} \int_M |S_{ij} + f_{ij} + \phi_{ij} + \frac{a+c}{2} g_{ij}|^2 e^{-f} d\mu.
\]

\[
-\frac{n(a+c)^2}{8} - \frac{1}{2} \int_M |\phi_{ij}|^2 e^{-f} d\mu.
\]

\[
-\int_M S_{ij} \phi_{ij} e^{-f} d\mu + \left(2(b-d) - \frac{1}{2}\right) \int_M |S_{ij}|^2 e^{-f} d\mu.
\]

\[
(2(b-d) - \frac{1}{2}) \int_M \alpha |\Delta \phi|^2 e^{-f} d\mu.
\]

\[
+\frac{1}{2} \int_M \alpha |\Delta \phi - \langle \nabla \phi, \nabla f \rangle|^2 e^{-f} d\mu.
\]

\[
+\int_M \alpha \langle \nabla \phi, \nabla f \rangle |\nabla \phi|^2 e^{-f} d\mu.
\]
\end{proof}
\subsection{Proof of Theorem 2}
\begin{proof}

From the Theorem 1, we have:

\[
\frac{d}{dt} \lambda^{b'}_{a'} = \frac{1}{2} \int_M \left| S_{ij} + f_{ij} + \phi_{ij} + \frac{a + c}{2} g_{ij} \right|^2 e^{-f} d\mu - \frac{n(a + c)^2}{8} - \frac{1}{2} \int_M |\phi_{ij}|^2 e^{-f} d\mu - \int_M S_{ij} \phi_{ij} e^{-f} d\mu
\]

\[
+ \left(2(b - d) - \frac{1}{2}\right) \int_M |S_{ij}|^2 e^{-f} \, d\mu 
+ \left(2(b - d) - \frac{1}{2}\right) \int_M \alpha |\Delta \phi|^2 e^{-f} \, d\mu 
\]
\[
+ \frac{1}{2} \int_M \alpha |\Delta \phi - \langle \nabla \phi, \nabla f \rangle|^2 e^{-f} \, d\mu 
+ \int_M \alpha \langle \nabla \phi, \nabla f \rangle |\nabla \phi|^2 e^{-f} \, d\mu.
\]

We want to analyze the condition under which the following quantity is nondecreasing:

\[
\lambda^{b'}_{a'}(t) + \frac{n(a + c)^2}{8} t.
\]

Compute the derivative with respect to \(t\):

\[
\frac{d}{dt} \left( \lambda^{b'}_{a'}(t) + \frac{n(a + c)^2}{8} t \right) = \frac{d}{dt} \lambda^{b'}_{a'} + \frac{n(a + c)^2}{8}.
\]

Substitute the expression for \(\frac{d}{dt} \lambda^{b'}_{a'}\):

\[
\frac{d}{dt} \left( \lambda^{b'}_{a'}(t) + \frac{n(a + c)^2}{8} t \right) = \frac{1}{2} \int_M \left| S_{ij} + f_{ij} + \phi_{ij} + \frac{a + c}{2} g_{ij} \right|^2 e^{-f} d\mu - \frac{n(a + c)^2}{8} + \frac{n(a + c)^2}{8}
\]

\[
- \frac{1}{2} \int_M |\phi_{ij}|^2 e^{-f} d\mu - \int_M S_{ij} \phi_{ij} e^{-f} d\mu + \left(2(b - d) - \frac{1}{2}\right) \int_M |S_{ij}|^2 e^{-f} d\mu
\]

\[
+ \left(2(b - d) - \frac{1}{2}\right) \int_M \alpha |\Delta \phi|^2 e^{-f} d\mu + \frac{1}{2} \int_M \alpha |\Delta \phi - \langle \nabla \phi, \nabla f \rangle|^2 e^{-f} d\mu + \int_M \alpha \langle \nabla \phi, \nabla f \rangle |\nabla \phi|^2 e^{-f} d\mu.
\]

Simplify the expression:

\[
\frac{d}{dt} \left( \lambda^{b'}_{a'}(t) + \frac{n(a + c)^2}{8} t \right) = \frac{1}{2} \int_M \left| S_{ij} + f_{ij} + \phi_{ij} + \frac{a + c}{2} g_{ij} \right|^2 e^{-f} d\mu - \frac{1}{2} \int_M |\phi_{ij}|^2 e^{-f} d\mu
\]

\[
- \int_M S_{ij} \phi_{ij} e^{-f} d\mu + \left(2(b - d) - \frac{1}{2}\right) \int_M |S_{ij}|^2 e^{-f} d\mu + \left(2(b - d) - \frac{1}{2}\right) \int_M \alpha |\Delta \phi|^2 e^{-f} d\mu
\]

\[
+ \frac{1}{2} \int_M \alpha |\Delta \phi - \langle \nabla \phi, \nabla f \rangle|^2 e^{-f} d\mu + \int_M \alpha \langle \nabla \phi, \nabla f \rangle |\nabla \phi|^2 e^{-f} d\mu.
\]

Using the condition:

\[
|S_{ij} - \frac{1}{4(b - d) - 1} \phi_{ij}| \geq 2 \sqrt{\frac{b - d}{4(b - d) - 1}} |\phi_{ij}|,
\]

we substitute to get:

\[
\frac{d}{dt} \left( \lambda^{b'}_{a'}(t) + \frac{n(a + c)^2}{8} t \right) \geq \frac{1}{2} \int_M \left( |S_{ij} - \frac{1}{4(b - d) - 1} \phi_{ij}|^2 - \frac{4(b - d)}{(4(b - d) - 1)^2} |\phi_{ij}|^2 \right) e^{-f} d\mu.
\]
\end{proof}

This shows that under the given condition, the quantity \(\lambda^{b'}_{a'}(t) + \frac{n(a + c)^2}{8} t\) is nondecreasing along the extended Ricci flow.
\subsection{Proof of Theorem 3}
\begin{proof}
We start by taking the modified equation under the normalized extended Ricci flow:

\[
-\Delta_\phi u + (a + c) u \log u + (b - d) S u = \lambda^{b'}_{a'}(g) u, \quad \text{with} \int_M u^2 d\mu = 1.
\]

Differentiate with respect to time \(t\):

\[
\frac{d}{dt} \lambda^{b'}_{a'} = \frac{d}{dt} \int_M \left( -u \Delta_\phi u + (a + c) u^2 \log u + (b - d) S u^2 \right) d\mu.
\]

Apply the Leibniz rule for differentiation:

\[
\frac{d}{dt} \lambda^{b'}_{a'} = \int_M \frac{\partial}{\partial t} \left( -\Delta_\phi u + (a + c) u \log u + (b - d) S u \right) u d\mu + \int_M \left( -\Delta_\phi u + (a + c) u \log u + (b - d) S u \right) \frac{\partial}{\partial t} (u d\mu).
\]

Using the derivatives under the normalized extended Ricci flow:

\[
\frac{\partial}{\partial t} \Delta_\phi u = 2 S_{ij} u_{ij} - 2 S_{ij} \phi_i u_j + \Delta_\phi u_t - (\phi_t)_i u^i - 2\alpha (\Delta \phi) \phi_i u^i - \frac{2}{n} r \Delta_\phi u.
\]

Then,
\[
\frac{d}{dt} \lambda^{b'}_{a'} = \int_M \Bigg\{ -\Bigg( 2 S_{ij} u_{ij} + \Delta_\phi u_t - 2 S_{ij} \phi_i u_j - (\phi_t)_i u^i
- 2 \alpha (\Delta \phi) \langle \nabla \phi, \nabla u \rangle - \frac{2}{n} r \Delta_\phi u \Bigg)
\]

\[
+(a+c) u_t + (a+c) u_t \log u + (b-d) S u_t + (b-d) S u_t \Bigg\} \, u \, d\mu.
\]

Simplify the expression:

\[
\frac{d}{dt} \lambda^{b'}_{a'} = \int_M \left\{ -2 S_{ij} u_{ij} u + 2 S_{ij} \phi_i u_j u + (\phi_t)_i u^i u + (b-d) S u^2 + (a+c) u u_t \right.
\]
\[
\left. + 2 \alpha (\Delta \phi) \langle \nabla \phi, \nabla u \rangle u + \frac{2}{n} r u \Delta_\phi u + (-\Delta_\phi u + (a+c) u \log u + (b-d) S u) u_t \right\} d\mu.
\]

Differentiating \(\int_M u^2 d\mu = 1\) with respect to \(t\), we have:

\[
\int_M u u_t d\mu = -\frac{1}{2} \int_M (\phi_t + S - r) u^2 d\mu.
\]

Substitute into the expression:

\[
\frac{d}{dt} \lambda^{b'}_{a'} = - \int_M 2 S_{ij} u_{ij} u d\mu + 2 \int_M S_{ij} \phi_i u_j u d\mu + \int_M (\phi_t)_i u^i u d\mu + \int_M (b-d) S u^2 d\mu
\]

\[
+ \frac{a+c}{2} \int_M (\phi_t + S - r) u^2 d\mu + 2 \alpha \int_M (\Delta \phi) \langle \nabla \phi, \nabla u \rangle u d\mu + 2(b-d) \alpha \int_M |\Delta \phi|^2 u^2 d\mu - \frac{2}{n} r \int_M (-\Delta_\phi u + (b - d) S u) u d\mu.
\]

After transforming \(u^2 = e^{-f}\) and using integration by parts:

\[
\frac{d}{dt} \lambda^{b'}_{a'} = \frac{1}{2} \int_M \left| S_{ij} + f_{ij} + \phi_{ij} + \frac{a+c}{2} g_{ij} \right|^2 e^{-f} d\mu - \frac{n(a+c)^2}{8} - \frac{1}{2} \int_M |\phi_{ij}|^2 e^{-f} d\mu
\]

\[
- 2\frac{r}{n} \lambda^{b'} - \frac{a+c}{2}r - \int_M S_{ij} \phi_{ij} e^{-f} d\mu + \left(2(b - d) - \frac{1}{2}\right) \int_M |S_{ij}|^2 e^{-f} d\mu
\]

\[
+ \left(2(b - d) - \frac{1}{2}\right) \int_M \alpha |\Delta \phi|^2 e^{-f} d\mu + \frac{1}{2} \int_M \alpha |\Delta \phi - \langle \nabla \phi, \nabla f \rangle|^2 e^{-f} d\mu + \int_M \alpha \langle \nabla \phi, \nabla f \rangle |\nabla \phi|^2 e^{-f} d\mu.
\]
\end{proof}
\subsection{Proof of Theorem 4}
\begin{proof}
From the  Theorem 3, we have:

\[
\frac{d}{dt} \lambda^{b'}_{a'} = \frac{1}{2} \int_M \left| S_{ij} + f_{ij} + \phi_{ij} + \frac{a+c}{2} g_{ij} \right|^2 e^{-f} d\mu - \frac{n(a+c)^2}{8} - \frac{1}{2} \int_M |\phi_{ij}|^2 e^{-f} d\mu - 2\frac{r}{n} \lambda^{b'} - \frac{a+c}{2}r
\]

\[
- \int_M S_{ij} \phi_{ij} e^{-f} d\mu + \left(2(b - d) - \frac{1}{2}\right) \int_M |S_{ij}|^2 e^{-f} d\mu + \left(2(b - d) - \frac{1}{2}\right) \int_M \alpha |\Delta \phi|^2 e^{-f} d\mu
\]

\[
+ \frac{1}{2} \int_M \alpha |\Delta \phi - \langle \nabla \phi, \nabla f \rangle|^2 e^{-f} d\mu + \int_M \alpha \langle \nabla \phi, \nabla f \rangle |\nabla \phi|^2 e^{-f} d\mu.
\]

We want to analyze the condition under which the following quantity is increasing:

\[
\lambda^{b'}_{a'}(t) + \frac{n(a+c)^2}{8} t + \int_0^t \left( \frac{a+c}{2} r + \frac{2r}{n} \lambda^{b'} \right) (\tau) \, d\tau.
\]

Compute the derivative with respect to \(t\):

\[
\frac{d}{dt} \left( \lambda^{b'}_{a'}(t) + \frac{n(a+c)^2}{8} t + \int_0^t \left( \frac{a+c}{2} r + \frac{2r}{n} \lambda^{b'} \right) (\tau) \, d\tau \right).
\]

This expands to:

\[
\frac{d}{dt} \lambda^{b'}_{a'}(t) + \frac{n(a+c)^2}{8} + \left( \frac{a+c}{2} r + \frac{2r}{n} \lambda^{b'} \right).
\]

Substitute the expression for \(\frac{d}{dt} \lambda^{b'}_{a'}\) from Theorem 3.1':

\[
\frac{d}{dt} \left( \lambda^{b'}_{a'}(t) + \frac{n(a+c)^2}{8} t + \int_0^t \left( \frac{a+c}{2} r + \frac{2r}{n} \lambda^{b'} \right) (\tau) \, d\tau \right) = \frac{1}{2} \int_M \left| S_{ij} + f_{ij} + \phi_{ij} + \frac{a+c}{2} g_{ij} \right|^2 e^{-f} d\mu
\]

\[
- \frac{n(a+c)^2}{8} + \frac{n(a+c)^2}{8} - \frac{1}{2} \int_M |\phi_{ij}|^2 e^{-f} d\mu - \int_M S_{ij} \phi_{ij} e^{-f} d\mu + \left(2(b - d) - \frac{1}{2}\right) \int_M |S_{ij}|^2 e^{-f} d\mu
\]

\[
+ \left(2(b - d) - \frac{1}{2}\right) \int_M \alpha |\Delta \phi|^2 e^{-f} d\mu + \frac{1}{2} \int_M \alpha |\Delta \phi - \langle \nabla \phi, \nabla f \rangle|^2 e^{-f} d\mu + \int_M \alpha \langle \nabla \phi, \nabla f \rangle |\nabla \phi|^2 e^{-f} d\mu
\]

\[
- 2\frac{r}{n} \lambda^{b'} - \frac{a+c}{2}r + \frac{a+c}{2}r + \frac{2r}{n} \lambda^{b'}.
\]

Simplify the expression:

\[
\frac{d}{dt} \left( \lambda^{b'}_{a'}(t) + \frac{n(a+c)^2}{8} t + \int_0^t \left( \frac{a+c}{2} r + \frac{2r}{n} \lambda^{b'} \right) (\tau) \, d\tau \right) = \frac{1}{2} \int_M \left| S_{ij} + f_{ij} + \phi_{ij} + \frac{a+c}{2} g_{ij} \right|^2 e^{-f} d\mu
\]

\[
- \frac{1}{2} \int_M |\phi_{ij}|^2 e^{-f} d\mu - \int_M S_{ij} \phi_{ij} e^{-f} d\mu + \left(2(b - d) - \frac{1}{2}\right) \int_M |S_{ij}|^2 e^{-f} d\mu
\]

\[
+ \left(2(b - d) - \frac{1}{2}\right) \int_M \alpha |\Delta \phi|^2 e^{-f} d\mu + \frac{1}{2} \int_M \alpha |\Delta \phi - \langle \nabla \phi, \nabla f \rangle|^2 e^{-f} d\mu + \int_M \alpha \langle \nabla \phi, \nabla f \rangle |\nabla \phi|^2 e^{-f} d\mu.
\]

Using the condition:

\[
|S_{ij} - \frac{1}{4(b - d) - 1} \phi_{ij}| \geq 2 \sqrt{\frac{b - d}{4(b - d) - 1}} |\phi_{ij}|,
\]

we substitute to get:
\[
\frac{d}{dt} \left( \lambda^{b'}_{a'}(t) + \frac{n(a+c)^2}{8} t + \int_0^t \left( \frac{a+c}{2} r + \frac{2r}{n} \lambda^{b'} \right) (\tau) \, d\tau \right) \geq 
\]

\[
\frac{1}{2} \int_M \left( \left|S_{ij} - \frac{1}{4(b - d) - 1} \phi_{ij}\right|^2 - \frac{4(b - d)}{(4(b - d) - 1)^2} |\phi_{ij}|^2 \right) e^{-f} \, d\mu.
\]
This shows that under the given condition, the quantity

\[
\lambda^{b'}_{a'}(t) + \frac{n(a+c)^2}{8} t + \int_0^t \left( \frac{a+c}{2} r + \frac{2r}{n} \lambda^{b'} \right) (\tau) \, d\tau
\]

is increasing along the normalized extended Ricci flow.
\end{proof}

\end{document}